\documentclass{article}
\usepackage{amsfonts}
\begin{document}
\title { Another proofs of the geometrical forms of Paley-Wiener theorems for the
 Dunkl transform and inversion formulas for the Dunkl
intertwining operator and for its dual}
\date{}
\author{Khalifa TRIM\`ECHE\\
{\small Faculty of Sciences of Tunis, Department of
Mathematics},\\ {\small CAMPUS,  1060 Tunis, Tunisia}\\ email :
{\small khlifa.trimeche@fst.rnu.tn}} \maketitle
\begin{abstract}
In this paper we present another proofs of the geometrical forms
of Paley-Wiener theorems for the Dunkl transform given in [15],
and we prove inversion formulas for the Dunkl intertwining
operator $V_k$ and for its dual ${}^tV_k$ and we deduce the
expression of the representing distributions of the inverse
operators $V^{-1}_k$ and ${}^tV_k^{-1}$.
\end{abstract}
\section{Introduction} We consider the differential-difference
operators $T_j, j = 1, 2,\cdots,d$, on $\mathbb{R}^d$ introduced
by C.F.Dunkl in [3]. These operators are very important  in pure
mathematics and in Physics. They provide a useful tool in the
study of special functions with root systems [4,6, 2].  Moreover
the commutative algebra generated by these operators has been used
in the study of certain exactly solvable models of quantum
mechanics, namely the Calogero-Sutherland-Moser models, which deal
with systems of identical particles in a one dimensional space
(see [8,11, 12]).

C.F.Dunkl has proved in [5] that there exists a unique isomorphism
$V_k$ from the space of homogeneous polynomials $\mathcal{P}_n$ on
$\mathbb{R}^d$ of degree $n$ onto itself satisfying the
permutation relations $$T_jV_k = V_k \frac{\partial}{\partial
x_j},\quad j = 1, 2,\cdots,d,\eqno{(1.1)}$$ and $$V_k(1) =
1.\eqno{(1.2)}$$ This operator is called Dunkl intertwining
operator. It has been extended to an isomorphism from
$\mathcal{E}(\mathbb{R}^d)$ (the space of $C^\infty$-functions on
$\mathbb{R}^d)$ onto itself satisfying the relations (1.1) and
(1.2) (see[14]).

The operator $V_k$ possesses the integral  representation
$$\forall\; x \in \mathbb{R}^d,\quad V_k(f)(x) =
\int_{\mathbb{R}^d}f(y)d\mu_x(y), \quad f \in
\mathcal{E}(\mathbb{R}^d),\eqno{(1.3)}$$ where $\mu_x$ is a
probability measure on $\mathbb{R}^d$ with support in the closed
ball $B(0, \|x\|)$ of center $0$ and radius $\|x\|$ (see [13][14])

We have shown in [14] that for each $x \in \mathbb{R}^d$, there
exists a unique distribution $\eta_x$ in
$\mathcal{E}'(\mathbb{R}^d)$ (the space of distributions on
$\mathbb{R}^d$ of compact support) with support in $B(0, \|x\|)$
such that $$V_k^{-1}(f)(x) = \langle \eta_x,f\rangle, f \in
\mathcal{E}(\mathbb{R}^d)$$

We have studied also in [14] the transposed operator ${}^tV_k$ of
the operator $V_k$. It has the integral representation $$\forall\;
y \in \mathbb{R}^d,\;\; {}^tV_k(f)(y) =
\int_{\mathbb{R}^d}f(x)d\nu_y(x).\eqno{(1.4)} $$ where $\nu_y$ is
a positive measure on $\mathbb{R}^d$ with support in the set $\{x
\in \mathbb{R}^d /\;\|x\| \geq \|y\|\}$ and $f$ in
$\mathcal{D}(\mathbb{R}^d)$ (the space of $C^\infty$-functions on
$\mathbb{R}^d$ with compact support).

This operator is called Dual Dunkl intertwining operator.

We have proved in [14] that the operator ${}^tV_k$ is an
isomorphism from $\mathcal{D}(\mathbb{R}^d)$ onto itself,
satisfying the transmutation relations $$\forall\; y \in
\mathbb{R}^d,\;\;\; {}^tV_k(T_jf)(y) = \frac{\partial}{\partial
y_j} {}^tV_k(f)(y),\quad j = 1, 2,\cdots,d.\eqno{(1.5)}$$

Using the operator $V_k$ C.F.Dunkl has defined in [5] the Dunkl
kernel $K$ by $$\forall\; x \in \mathbb{R}^d,\; \forall\; z \in
\mathbb{C}^d,\; K(x, -iz) = V_k(e^{-i\langle
.,z\rangle})(x).\eqno{(1.6)}$$ Using this kernel C.F.Dunkl has
introduced in [5] a Fourier transform $\mathcal{F}_D$ called Dunkl
transform.

In this paper we proesent another proofs of the geometric forms of
Paley-Wiener theorems for the transform $\mathcal{F}_D$, given in
[15] p. 32-33, and we establish the following inversion formulas
for the operators $V^{-1}_k$ and ${}^tV^{-1}_k$ : $$\forall\;\; x
\in \mathbb{R}^d,\;\; V^{-1}_k(f)(x) = {}^tV_k(Q(f)) (x),\;\; f
\in \mathcal{D}(\mathbb{R}^d),\eqno{(1.7)} $$ $$\forall\; x \in
\mathbb{R}^d,\;\; {}^tV^{-1}_k(f)(x) = V_k(P(f))(x),\quad f \in
\mathcal{D}(\mathbb{R}^d),\eqno{(1.8)} $$ where $P$ and $Q$ are
pseudo-differential operators on $\mathbb{R}^d$.\\ From these
relations we deduce the expression of the representing
distributions $\eta_x$ and $Z_x$ of the inverse operators
$V^{-1}_k$ and ${}^tV^{-1}_k$ by using the representing measures
$\mu_x$ and $\nu_x$ of $V_k$ and ${}^tV_k$. They are given  by the
following formulas $$\forall\; x \in \mathbb{R}^d,\quad \eta_x =
{}^tQ(\nu_x),\eqno{(1.9)}$$ $$\forall\; x \in \mathbb{R}^d,\quad
Z_x = {}^tP(\mu_x),\eqno{(1.10)}$$ where ${}^tP$ and ${}^tQ$ are
the transposed of the oprators $P$ and $Q$.

The contents of the paper are as follow.

In section two we recall some basic facts from Dunkl's theory, we
discribe Dunkl operators and the Dunkl kernel.

We introduce in the third section  the Dunkl intertwining operator
$V_k$ and its dual ${}^tV_k$  and we present their properties.

We define int the fourth section the Dunkl  transform introduced
in [5] by C.F.Dunkl, and we give the main theorems proved for this
transform.

In the fifth section we give another proofs of the geometrical
forms of Paley-Wiener theorems for the Dunkl transform. The first
proofs of these theorems have been given in [15] p. 32-33. Next we
present some applications of the first theorem.

The sixth section is devoted to prove inversion formulas for the
Dunkl intertwining operator $V_k$ and for its dual ${}^tV_k$, and
we deduce the expression of the representing distributions of the
inverse operators $V^{-1}_k$ and ${}^tV^{-1}_k$.
\section{The eigenfunction of the Dunkl operators}

In this section we collect some notations and results on Dunkl
operators and the Dunkl kernel (see [4,5, 7, 9, 10]).
\subsection{Reflection Groups, Root Systems and Multiplicity
Functions}

We consider $\mathbb{R}^d$ with the euclidean scalar product
$\langle .,.\rangle$ and $\|x\| = \sqrt{\langle x , x\rangle}$. On
$\mathbb{C}^d, \|.\|$ denotes also the standard Hermitian norm,
while $\langle z, w\rangle = \sum^d_{j=1}z_j \overline{w_j}$ .

For $\alpha \in \mathbb{R}^d \backslash \{0\}$, let
$\sigma_\alpha$ be the reflection in the hyperplan $H_\alpha
\subset \mathbb{R}^d$ orthogonal to $\alpha$, i.e.
$$\sigma_\alpha(x) = x - \left(\frac{2\langle
\alpha,x\rangle}{\|\alpha\|^2}\right)\alpha .\eqno{(2.1)}$$ A
finite set $R \subset \mathbb{R}^d \backslash \{0\}$ is called a
root  system if $R \cap \mathbb{R}\alpha = \{\pm \alpha\}$ and
$\sigma_\alpha R = R$ for all $\alpha \in R$. For a give root
system $R$ the reflections $\alpha_\alpha, \alpha \in R$, generate
a finite group $W \subset O(d)$, the reflection group associated
with $R$. All reflections in $W$ correspond to suitable pairs of
roots. For a given $\beta \in \mathbb{R}^d \backslash \cup_{\alpha
\in R }H_\alpha$, we fix the positive subsystem $R_+ = \{\alpha
\in R ; \langle \alpha, \beta\rangle > 0\}$, then for each $\alpha
\in R$ either $\alpha \in R_+$ or $-\alpha \in R_+$.

A function $k : R \rightarrow \mathbb{C}$ on a root system $R$ is
called a multiplicity function if it is invariant under the action
of the associated reflection group $W$. If one regards $k$ as a
function on the corresponding reflections, this means that $k$ is
constant on the conjugacy classes of reflections in $W$. For
abbreviation, we introduce the index $$\gamma = \gamma(R) =
\sum_{\alpha \in R_+} = \sum_{\alpha \in R_+}
k(\alpha).\eqno{(2.2)}$$ Moreover, let $\omega_k$ denotes the
weight function $$\omega_k(x) = \prod_{\alpha \in R_+}|\langle
\alpha,x\rangle|^{2k(\alpha)}.\eqno{(2.3)}$$ which is
$W$-invariant and homogeneous of degree $2\gamma$.

For $d = 1$ and $W= \mathbb{Z}_2$, the multiplicity function $k$
is a single parameter denoted $\gamma > 0$ and $$\forall\; x \in
\mathbb{R},\;\; \omega_k(x) = |x|^{2\gamma}.\eqno{(2.4)} $$ We
introduce the Mehta-type constant $$c_k = \left(
\int_{\mathbb{R}^d}
e^{-\|x\|^2}\omega_k(x)dx\right)^{-1}.\eqno{(2.5)}$$ which is
known for all Coxeter groups $W$ (see [3, 6])
\subsection{Dunkl Operators and Dunkl kernel}

The Dunkl operators $T_j, j = 1,\cdots,d$, on $\mathbb{R}^d$,
associated with the finite reflection group $W$ and multiplicity
function $k$, are given for a function $f$ of class $C^1$ on
$\mathbb{R}^d$ by $$T_jf(x) = \frac{\partial}{\partial x_j}f(x) +
\sum_{x \in R_+} k(\alpha)\alpha_j \frac{f(x) -
f(\sigma_\alpha(x))}{\langle\alpha,x\rangle}.\eqno{(2.6)} $$ In
the case $k = 0$, the $T_j, j = 1, 2,\cdots,d$, reduce to the
corresponding partial derivatives. In this paper, we will assume
throughout that $k \geq 0$ and $\gamma > 0$.

For $f$ of class $C^1$ on $\mathbb{R}^d$ with compact support and
$g$ of class $C^1$ on $\mathbb{R}^d$ we have
$$\int_{\mathbb{R}^d}T_jf(x) g(x) \omega_k(x)dx =
-\int_{\mathbb{R}^d} f(x) T_j g(x) \omega_k(x) dx,\quad j = 1,
2,\cdots, d.\eqno{(2.7)}$$ For $y \in \mathbb{R}^d$, the system
$$\left\{ \begin{array}{ll} T_j u(x, y) &= y_j u(x, y),\quad j =
1, 2,\cdots,d,\\ u(0,y) &= 1,
\end{array}\right.\eqno{(2.8)}$$
admits a unique analytic solution on $\mathbb{R}^d$, denoted by
$K(x,y)$ and called Dunkl kernel.

This kernel has a unique holomorphic extension to $\mathbb{C}^d
\times \mathbb{C}^d$.\\ {\bf Example}

If $d = 1$ and $W = \mathbb{Z}_2$, the Dunkl kernel is given by
$$K(z,t) = j_{\gamma-1/2} (izt) + \frac{zt}{2\gamma+1}
j_{\gamma+1/2} (izt),\quad z, t \in \mathbb{C},\eqno{(2.9)}$$
where for $\alpha `\geq - 1/2, j_\alpha$ is the normalized Bessel
function defined by $$j_\alpha(u) = 2^\alpha \Gamma(\alpha + 1)
\frac{J_\alpha(u)}{u^\alpha} = \Gamma(\alpha +1) \sum^\infty_{n=0}
\frac{(-1)^n(u/2)^{2n}}{n!\Gamma(n+\alpha +1)}, \quad u \in
\mathbb{C},\eqno{(2.10)} $$ with $J_\alpha$ is the Bessel function
of first kind and index $\alpha$ (see [5]).

The Dunkl kernel possesses the following properties.
\begin{itemize}
\item[(i)] For $z, t \in \mathbb{C}^d$, we have $K(z,t) = K(t,z), K(z,0) =
1$, and $K(\lambda z,t) = K(z, \lambda t)$ for all $\lambda \in
\mathbb{C}$.
\item[(ii)] For all $\nu \in \mathbb{Z}^d_+, x \in \mathbb{R}^d$,
and $z \in \mathbb{C}^d$ we have $$|D^\nu_z K(x,z)| \leq
\|x\|^{|\nu|} \exp \left[\max_{w \in W} \langle wx, Re
z\rangle\right].\eqno{(2.11)}$$ In particular $$|D^\nu_z K(x,z)|
\leq \|x\|^{|\nu|} \exp [\|x\| \|Re z\|]],\eqno{(2.12)} $$
$$|K(x,z)| \leq \exp [\|x\| \|Re z\|]\eqno{(2.13)}$$ and for all
$x, y \in \mathbb{R}^d$ : $$|K(ix,y)| \leq 1, \eqno{(2.14)}$$ with
$$D^\nu_z = \frac{\partial^{|\nu|}} {\partial^{\nu_1}z_1 \cdots
\partial z^{\nu_d}_d } \mbox{ and } |\nu| = \nu_1 + \cdots+
\cdots + \nu_d.$$ \item[(iii)] For all $x, y \in \mathbb{R}^d$
 and $w \in W$ we have
 $$K(-ix, y) = \overline{K(ix, y)} \mbox{ and } K(wx, wy) =
  K(x,y).\eqno{(2.15)}$$
 \item[(iv)] The function $K(x,z)$ admits for all $x \in
 \mathbb{R}^d$ and $z \in \mathbb{C}^d$ the following Laplace type
 integral representation
 $$K(x,z) = \int_{\mathbb{R}^d} e^{\langle y,z\rangle}d\mu_x(y),
 \eqno{(2.16)}$$
 where $\mu_x$ is a probability measure on $\mathbb{R}^d$ with
 support in the closed ball $B(0, \|x\|)$ of center $0$ and radius
 $\|x\|$. Moreover we have
 $$Supp \mu_x \cap \{y \in \mathbb{R}^d / \|y\| = \|x\|\} \neq
 \emptyset,\eqno{(2.17)}$$
 (see [13]).
\end{itemize}
{\bf Remark 2.1}

When $d = 1$ and $W = \mathbb{Z}_2$, the relation (2.16) is of the
form $$K(x,z) = \frac{\Gamma(\gamma+1/2)}
{\sqrt{\pi}\Gamma(\gamma)}|x|^{-2\gamma}\int^{|x|}_{-|x|}(|x| -
y)^{\gamma-1}(|x| + y)^\gamma e^{yz}dy.\eqno{(2.18)}$$ Then in
this case the measure $\mu_x$ is given for all $x \in
\mathbb{R}\backslash \{0\}$ by $d\mu_x(y) = \kappa(x,y) dy$ with
$$\kappa(x,y) =
\frac{\Gamma(\gamma+1/2)}{\sqrt{\pi}\Gamma(\gamma)}|x|^{-2\gamma}(|x|
-y)^{\gamma-1}(|x|+y)^\gamma 1_{]-|x|, |x|[}(y),\eqno{(2.19)} $$
where $1_{]-|x|, |x|[}$ is the characteristic function of the
interval $]-|x|, |x|[$.

We remark that by change of variables,  the relation (2.18) takes
the following form $$\forall\; x \in \mathbb{R}^d, \; \forall\; z
\in \mathbb{C}^d,\;\; K(x,z) = \frac{\Gamma(\gamma + 1/2)}
{\sqrt{\pi}(\gamma)}\int^1_{-1}e^{txz}(1-t^2)^{\gamma-1}(1+t)dt,\eqno{(2.20)}
$$
\section{The Dunkl intertwining operator and its dual}
{\bf Notation} We denote by

- $C(\mathbb{R}^d)$ (resp. $C_c(\mathbb{R}^d))$ the space of
continuous functions on $\mathbb{R}^d$ (resp. with compact
support).

- $C^p(\mathbb{R}^d)$ (resp. $C^p_c(\mathbb{R}^d))$ the space of
functions of class $C^p$ on $\mathbb{R}^d$ (resp. with compact
support).

- $\mathcal{E}(\mathbb{R}^d)$ the space of $C^\infty$-functions on
$\mathbb{R}^d$.

- $\mathcal{D}(\mathbb{R}^d)$ the space of $C^\infty$-functions on
$\mathbb{R}^d$ with compact support.

- $\mathcal{S}(\mathbb{R}^d)$ the space of $C^\infty$-functions on
$\mathbb{R}^d$ which are rapidly decreasing as their derivatives.

We provide these spaces with the classical topology. We consider
also the following spaces.

- $\mathcal{E}'(\mathbb{R}^d)$ the space of distributions on
$\mathbb{R}^d$ with compact support ; it is the topological dual
of $\mathcal{E}(\mathbb{R}^d)$.

- $\mathcal{S}'(\mathbb{R}^d)$ the space of tempered distributions
on $\mathbb{R}^d$ ; it is the topological dual of
$\mathcal{S}(\mathbb{R}^d)$.\vspace{5mm}

The Dunkl intertwining operator $V_k$ is defined on
$\mathcal{C}(\mathbb{R}^d)$ by $$\forall\; x \in \mathbb{R}^d,
V_k(f)(x) = \int_{\mathbb{R}^d}f(y)d\mu_x(y),\eqno{(3.1)}$$ where
$\mu_x$ is the measure given by the relation (2.16) (see [14]).

We have $$\forall\; x \in \mathbb{R}^d,\;\; \forall\; z \in
\mathbb{C}^d,\;\; K(x,z) = V_k(e^{<.,y>})(x).\eqno{(3.2)}$$ The
operator ${}^tV_k$ satisfying for $f$ in $C_c(\mathbb{R}^d)$ and
$g$ in $C(\mathbb{R}^d)$, the relation
$$\int_{\mathbb{R}^d}{}^tV_k(f)(y)g(y)dy =
\int_{\mathbb{R}^d}V_k(g)(x) f(x) \omega_k(x)dx.\eqno{(3.3)}$$ is
given by $$\forall\; y \in \mathbb{R}^d, {}^tV_k(f)(y) =
\int_{\mathbb{R}^d}f(x) d\nu_y(x),\eqno{(3.4)}$$ where $\nu_y$ is
a positive measure on $\mathbb{R}^d$ whose support satisfies
$$Supp \nu_y \subset \{x \in \mathbb{R}^d / \|x\|\geq \|y\|\}
\mbox{ and } Supp \nu_y \cap \{x \in \mathbb{R}^d / \|x\| =
\|y\|\} \neq \emptyset. \eqno{(3.5)} $$ This operator is called
the dual Dunkl intertwining operator (see [14]). \vspace{3mm}

 The following theorems give some properties
of the operators $V_k$ and ${}^tV_k$ (see [14]).\newpage {\bf
Theorem 3.1}
\begin{itemize}
\item[(i)] The operator $V_k$ is a topological isomorphism from
$\mathcal{E}(\mathbb{R}^d)$ onto itself satisfying the
transmutation relations $$\forall\; x \in \mathbb{R}^d,\;\;
T_jV_k(f)(x) = V_k\left(\frac{\partial}{\partial
y_j}f\right)(x),\;\; j = 1, 2,\cdots,d,\;\; f \in
\mathcal{E}(\mathbb{R}^d).\eqno{(3.6)}$$
\item[(ii)] For each $x \in \mathbb{R}^d$, there exists a unique
distribution $\eta_x$ in $\mathcal{E}'(\mathbb{R}^d)$ with support
in the ball $B(0, \|x\|)$ such that for all $f$ in
$\mathcal{E}(\mathbb{R}^d)$ we have $$V^{-1}_k(f)(x) = \langle
\eta_x, f\rangle.\eqno{(3.7)}$$
\end{itemize}
Moreover $$Supp \eta_x \cap \{y \in \mathbb{R}^d / \|y\| = \|x\|\}
\neq \emptyset. \eqno{(3.8)} $$ {\bf Theorem 3.2}
\begin{itemize}
\item[(i)] The operator ${}^tV_k$ is a topological isomorphism
form $\mathcal{D}(\mathbb{R}^d)$ (resp.
$\mathcal{S}(\mathbb{R}^d))$ onto itself, satisfying the
transmutation relations $$\forall\; y \in \mathbb{R}^d,\;
{}^tV(T_jf)(y) = \frac{\partial}{\partial y_j}{}^tV(f)(y),\; j =
1, 2,\cdots,d, f \in \mathcal{D}(\mathbb{R}^d).\eqno{(3.9)} $$
\item[ii)] For each $y \in \mathbb{R}^d$, there exists a unique
distribution $Z_y$ in $\mathcal{S}'(\mathbb{R}^d)$ with support in
the set $\{x \in \mathbb{R}^d / \|x\| \geq \|y\|\}$ such that for
all $f$ in $\mathcal{D}(\mathbb{R}^d)$ we have
$${}^tV_k^{-1}(f)(y) = \langle Z_y, f\rangle.\eqno{(3.10)}$$
Moreover $$Supp Z_y \cap \{ x \in \mathbb{R} / \|x\| = \|y\|\}
\neq \emptyset.\eqno{(3.11)}$$
\end{itemize}
{\bf Example 3.1}

When $d = 1$ and $W = \mathbb{Z}_2$, the Dunkl intertwining
operator $V_k$ is defined by (3.1) with for all $x \in \mathbb{R}
\backslash \{0\}, d\mu_x(y) = \kappa(x,y) dy$, where $\kappa$
given by the relation (2.19).

The dual Dunkl intertwining operator ${}^tV_k$ is defined by (3.4)
with $d\nu_y(x) = \kappa(x,y) \omega_k(x)dx$, where $\kappa$ and
$\omega_k$ given respectively by the relations (2.19) and (2.4).\\
{\bf Example 3.2}

The Dunkl intertwining operator $V_t$ of index $\gamma =
\sum^d_{i=1} \alpha_i, \alpha_i > 0$, associated with the
reflection group $\mathbb{Z}_2 \times \mathbb{Z}_2 \times \cdots
\times \mathbb{Z}_2$ on $\mathbb{R}^d$, is given for all $f$ in
$\mathcal{E}(\mathbb{R}^d)$ and for all $x \in \mathbb{R}^d$ by
$$V_k(f)(x) = \prod^d_{i=1}\left(\frac{\Gamma(\alpha_i + 1/2)}
{\sqrt{\pi}\Gamma(\alpha_i)}\right) \int_{[-1,1]d}f(t_1x_1,
t_2x_2,\cdots,t_dx_d)$$ \hspace{4cm}$ \times
\displaystyle{\prod^d_{i=1}}(1-t^2_i)^{\alpha_i-1}(1+t_i)dt_1\cdots
dt_d,$\hfill(3.12)\\ (see [16])
\section{Dunkl transform}

In this section we define the Dunkl transform and we give the main
results satisfied by this transform (see [5, 9, 10]).\\ {\bf
Notation} We denote by

- $L^p_k(\mathbb{R}^d), p \in [1, + \infty]$, the space of
measurable functions on $\mathbb{R}^d$ such that
\begin{eqnarray*}
\|f\|_{k,p} &=& \left(\int_{\mathbb{R}^d}|f(x)|^p
\omega_k(x)dx\right)^{1/p}  < + \infty,\;\; \mbox{ if } 1\leq p <
+ \infty,\\ \|f\|_{k,\infty} &=& \displaystyle{ess\sup_{x \in
\mathbb{R}^d}}|f(x)| < + \infty.
\end{eqnarray*}

- $H(\mathbb{C}^d)$ the space of entire functions on
$\mathbb{C}^d$ which are rapidly decreasing and of exponential
type.\vspace{3mm}\\
 The Dunkl transform of a function $f$ in
$\mathcal{D}(\mathbb{R}^d)$ is given by $$\forall\; y \in
\mathbb{R}^d,\;\; \mathcal{F}_D(f)(y) =
\int_{\mathbb{R}^d}f(x)K(x, -iy)\omega_k(x)dx.\eqno{(4.1)}$$ This
transform has the following properties.
\begin{itemize}
\item[i)] For $f$ in $L^1_k(\mathbb{R}^d)$ we have $\|\mathcal{F}_D(f)\|_{k, \infty}
\leq \|f\|_{k,1}$.
\item[(ii)] Let $f$ be in $\mathcal{D}(\mathbb{R}^d)$. If $f^-(x) = f(-x)$
and $f_w(x) = f(wx)$ for $x \in \mathbb{R}^d$,  $w \in W$, then
for all $y \in \mathbb{R}^d$ we have $$\mathcal{F}_D(f^-)(y) =
\overline{\mathcal{F}_D(f)(y)} \mbox{ and } \mathcal{F}_D(f_w)(y)
= \mathcal{F}_D(f)(wy).\eqno{(4.2)}$$
\item[iii)] For all $f$ in $\mathcal{S}(\mathbb{R}^d)$ we have
$$\mathcal{F}_D(f) = \mathcal{F}\,o\,{}^tV_k(f),\eqno{(4.3)}$$
\end{itemize}
where $\mathcal{F}$ is the classical Fourier transform on
$\mathbb{R}^d$ given by $$\forall\; y \in \mathbb{R}^d,\;\;
\mathcal{F}(f)(y) = \int_{\mathbb{R}^d}f(x) e^{-i\langle
x,y\rangle}dx,\;\; f\in \mathcal{D}(\mathbb{R}^d),\eqno{(4.4)}$$
The following theorems are proved in [9, 10].\\ {\bf Theorem 4.1.}
The transform $\mathcal{F}_D$ is a topological isomorphism
\begin{itemize}
\item[i)] from $\mathcal{D}(\mathbb{R}^d)$ onto
$\mathbb{H}(\mathbb{C}^d)$,
\item[ii)] from $\mathcal{S}(\mathbb{R}^d)$ onto itself.
\end{itemize}
The inverse transform is given by
 $$\forall\; x \in \mathbb{R}^d,\;\; \mathcal{F}^{-1}_D(h)(x)
=
\frac{c^2_k}{2^{2\gamma+d}}\int_{\mathbb{R}^d}h(y)K(x,iy)
\omega_k(y)dy.\eqno{(4.5)}$$\newpage \noindent{\bf Remark 4.1}

An other proof of Theorem 4.1 is given in [15].\\ {\bf Theorem
4.2.} Let $f$ be in $L^1_k(\mathbb{R}^d)$ such that the function
$\mathcal{F}_D(f)$ belongs to $L^1_k(\mathbb{R}^d)$. Then we have
the following intersion formula for the transform $\mathcal{F}_D$
: $$f(x) = \frac{c^2_k}{2^{2\gamma+d}}\int_{\mathbb{R}^d}
\mathcal{F}_D(f)(y)K(x, iy)\omega_k(y)dy,\;\; a.e.\eqno{(4.6)}$$
{\bf Theorem 4.3.}
\begin{itemize}
\item[i)] \underline{Plancherel formula for} $\mathcal{F}_D$. For all $f$ in
$\mathcal{D}(\mathbb{R}^d)$ we have $$\int_{\mathbb{R}^d}|f(x)|^2
\omega_x(x)dx =
\frac{c^2_k}{2^{2\gamma+d}}\int_{\mathbb{R}^d}|\mathcal{F}_D(f)(y)|^2
\omega_k(y)dy.\eqno{(4.7)} $$
\item[ii)] \underline{Plancherel Theorem for} $\mathcal{F}_D$. The
renormalized Dunkl transform $f\rightarrow 2^{-\gamma-d/2} \times
c_k \mathcal{F}_D(f)$ can be uniquely extended to an isometric
isomorphism on $L^2_k(\mathbb{R}^d)$.
\end{itemize}
\section{Another proofs of the geometrical
 forms of the Paley-Wiener theorems for the Dunkl
 transform}
In this section we present another proofs of the geometrical forms
of  Paley-Wiener theorems for the transform $\mathcal{F}_D$, given
in [15] p. 32-33.\vspace{5mm}

We define the indicator function $I_E$ of a compact subset $E$ of
$\mathbb{R}^d$ by $$\forall\; y \in \mathbb{R}^d, I_E(y) = \sup_{x
\in E }\langle x, y\rangle,$$ for example, if $E$ is the ball with
center $0$ and radius $R$, we have $$\forall\; y \in \mathbb{R}^d,
I_E(y) = R\|y\|. $$

The indicator of $E$ determines, for each hyperplane direction,
the smallest closed half space which contain $E$. Consequently the
Hahn-Banach theorem shows that $I_E$ determines the convex
envelope $\hat{E}$ of $E$ and $I_{\hat{E}} = I_E$.

In [15] p. 29 we find other properties of the indicator function
$I_E$.
\subsection{The Dunkl transform of functions}
{\bf Theorem 5.1.} Let $E$ be a $W$-invariant compact convex set
of $\mathbb{R}^d$ and $f$ an entire function on $\mathbb{C}^d$.
Then $f$ is the Dunkl transform of a function in
$\mathcal{D}(\mathbb{R}^d)$ with support in $E$, if and only if
for all $q \in \mathbb{N}$ there exists a positive constant $C_q$
such that $$\forall\; z \in \mathbb{C}^d, |f(z)| \leq C_q(1 +
\|z\|)^{-q}e^{I_E(Im z)}\eqno{(5.1)}$$ {\bf Proof}

The necessity condition of the theorem is proved in [9] corollary
4.10, p. 156.

We prove now the sufficiency condition.

The upper bounds (5.1) imply that the function $\varphi$ given by
$$\forall\; x \in \mathbb{R}^d,\; \varphi(x) =
\int_{\mathbb{R}^d}f(y)K(iy, x) \omega_k (y)dy, \eqno{(5.2)}$$ is
a $C^\infty$-function on $\mathbb{R}^d$.

We consider the function $$\forall\; x \in \mathbb{R}^d, \phi(x) =
\int_{\mathbb{R}^d}f(y) K(iy, x)
\frac{\omega_k(y)}{(1+\|y\|^2)^p}dy,\eqno{(5.3)}$$ with $p \in
\mathbb{N}$ such that $p > \gamma + \frac{d}{2} +1$.\\ This
function is of class $C^\infty$ on $\mathbb{R}^d$ and we have
$$\forall\; x \in \mathbb{R}^d,\; (I - \Delta_k)^p \phi(x) =
\varphi(x),\eqno{(5.4)}$$ where $\Delta_k = \sum^d_{j=1}T^2_j$ is
the Dunkl Laplacian. As $E$ is $W$-invariant then to show that
Supp$\varphi \subset E$, it is sufficient to prove that the
support of $\phi$ is contained in $E$. The relation (5.3) can also
be written in the form $$\forall\; x \in \mathbb{R}^d,\; \phi(x) =
\int_{\mathbb{R}^d}f(y)K(iy, x) M_k(y)dy,\eqno{(5.5)}$$ where
$M_k(y) = \frac{\omega_k(y)}{(1+\|y\|^2)^p}$.\\ From (2.3) we have
$$\int_{\mathbb{R}^d}M_k(y)dy < + \infty.\eqno{(5.6)}$$ We put
$$m_k(z) = m_k(z_1,\cdots,z_d) = \int_{\mathbb{R}^d}
\frac{M_k(y_1,\cdots,y_d)}{(y_1-z_1)\cdots(y_d-z_d)} dy_1\cdots
dy_d,\eqno{(5.7)}$$ and for $a_j \geq 0, j = 1, 2,\cdots,d$,
$$m_k^{a_1,\cdots,a_d}(z) =
\int^{a_1}_{-a_1}\cdots\int^{a_d}_{-a_d}\frac{M_k(y_1,\cdots,y_d)}{(y_1-z_1)\cdots(y_d-z_d)}
dy_1\cdots,dy_d.\eqno{(5.8)}$$ These functions are respectively
holomorphic on $\mathbb{C}^d \backslash \mathbb{R}^d$ and
$\mathbb{C}^d \backslash (\prod^d_{j=1}[-a_j, a_j])$, and for all
$z \in \mathbb{C}^d \backslash \mathbb{R}^d$ we have
$$\lim_{a_1,\cdots,a_2\rightarrow + \infty}m_k^{a_1,\cdots,a_2}(z)
= m_k(z).\eqno{(5.9)} $$ Using Riemann sums, the function
$m_k^{a_1,\cdots,a_d}$ can also be written in the form
\newpage
$m_k^{a_1\cdots,a_2}(z) = \displaystyle{\lim_{n_1,\cdots,n_d
\rightarrow + \infty}\left(\prod^d_{j=1} \frac{2a_j}{n_j}\right)}
\times $\\

 \hspace{3cm}$\displaystyle{\sum^{n_1}_{r_1 = 0} \cdots
\sum^{n_d}_{r_d =0} \frac{M_k(-a_1 + r_1
\frac{2a_1}{n_1},\cdots,-a_d + r_d
\frac{2a_d}{n_d})}{\prod^d_{j=1}(-a_j + r_j
\frac{2a_j}{n_j}-R_j)}}.$\hfill (5.10)\\ For $j = 1, 2,\cdots,d$,
let $\Gamma_j$ the quasi-rectangular path in $\mathbb{C}$
determined by the points $z_j = - R_j - ib_j, z_j = - R_j +
i\eta_j$, $z_j = + R_j + i\eta_j, z_j = R_j + ib_j$, with $R_j >
0, 0 < b_j < \frac{a_j}{n_j}\;, \eta_j > 0$, the half circles
$C_{rj}, r_j = 1, 2,\cdots,n_j$, of center $(-a_j + r_j
\frac{2a_j}{n_j})$
 and radius $\frac{a_j}{n_j}$ and the segments $\{z_j = x_j
 + ib_j, x_j \in [-R_j, -a_j(1 +
 \frac{1}{n_j})+b_j]\}$ $\{z_j = x_j + ib_j, x_j \in [a_j(1 +
 \frac{1}{n_j})- b_j,R_j]\}$.\\
 From Cauchy theorem we have
 $$I_{a_1, \cdots,a_d} = \int_{\Gamma_1}\cdots \int_{\Gamma_d}f(z)K(iz, x) m_k^{a_1,\cdots,a_d}
 (z)dz = 0.\eqno{(5.11)}$$
 On the other hand from (5.8) and (5.10) we have
 \begin{eqnarray*}
 I_{(a_1\cdots,a_d)} &=&\int^{a_1}_{-a_1}\cdots
 \int^{a_d}_{-a_d}M_k(y_1,\cdots,y_d)\int_{\Gamma_1}\cdots\int_{\Gamma_d}\frac{f(z)K(iz,x)}
 {\prod^d_{j=1}(y_j-z_j)}dz dy_1,\cdots, dy_d,\\
 &=& \renewcommand{\arraystretch}{0.5}
 \begin{array}[t]{c}\\
 lim\\
 {\scriptstyle n_1 \rightarrow + \infty}\\
 {\scriptstyle n_d \rightarrow + \infty}
 \end{array}
 \renewcommand{\arraystretch}{} \Big(\prod^d_{j=1}\frac{2a_j}
 {n_j}\Big)\sum^{n_1}_{r_1 =0}\cdots \sum^{n_d}_{r_d=0} M_k (-a_1
 + r_1 \frac{2a_1}{n_1},\cdots, - a_d + r_d \frac{2a_d}{n_d})\\
 &&\times \int_{\Gamma_1}\cdots\int_{\Gamma_d} \frac{f(z_1,\cdots,z_d)K(iz_1,\cdots,iz_d,x)}
 {\prod^d_{j=1}(-a_j + r_j \frac{2a_j}{n_j}-z_j)}
 dz_1\cdots dz_d .\;\;\;\;\;\;\;\;\;\;\;\; (5.12)
 \end{eqnarray*}
 We apply residus theorem to the integrals of the second member of
 this relation. As from (2.13) and (5.1) the integrals on the segments
 $\{z_j = - R_j + it_j, t_j \in [b_j, \eta_j]\}$ and
 $\{z_j = R_j + it_j, t_j \in [b_j, R_j]\}, j = 1, 2,\cdots,d$,
  tend to zero when
 $R_1,\cdots,R_d \rightarrow + \infty$. Then if we make
  $R_1, R_2,\cdots, R_d {\rightarrow +
 \infty}$, we deduce from (5.11) that
 $$\int^{-a_1-b_1}_{-\infty}\cdots\int^{-a_d-b_d}_{-\infty}
 f(y+ib)K(i(y+ib),x)
 m^{a_1,\cdots,a_1} (y+ib)dy +$$
 $$\int^\infty_{a_1 +b_1 }\cdots \int^{+ \infty}_{a_d+b_d}
 f( y +ib)K(i(y +ib),x) m^{a_1,\cdots,a_d}(y+ib)dy +$$
 $$\int^{a_1}_{-a_1}\cdots \int^{a_d}_{-a_d}f(y)K(iy, x)
 \frac{\omega_k(y)}{(1+\|y\|^2)^p}dy \;\;\;-$$
 $$\frac{1}{(i\pi)^d}\int_{\mathbb{R}^d}f(y + i\eta)
 K(i(y+i\eta),x)
 m^{a_1\cdots,a_d}(y + i\eta)dy = 0,
 \eqno{(5.13)}$$
 where $b = (b_1, b_2,\cdots,b_d)$ , $\eta = (\eta_1, \eta_2,\cdots,\eta_d)$and $y = (y_1, y_2,\cdots,y_d)$.
 But from (5.1) (2.13) and (5.6) we deduce that the two first
 integrals tend to zero when $a_1,\cdots,a_d \rightarrow +
 \infty$. Thus by making $a_1,\cdots,a_d \rightarrow + \infty$ in
 (5.13) we obtain
 $$\forall\; x \in \mathbb{R}^d, \phi(x) = \frac{1}{(i\pi)^d}\int_{\mathbb{R}}
 f(y + i\eta)K(i(y + i\eta),x) m_k(y + i\eta)dy.$$
 Using (5.1) and (2.13) we deduce that there exists a positive
 constant $M$ such that
 $$\forall\; x \in \mathbb{R}^d, |\phi(x)| \leq e^{I_E(\eta)-\langle x, \eta\rangle} M \int_{\mathbb{R}^d}
 \frac{|m_k(y + i\eta)|}{(1+\|y\|^2)^p}dy, \eqno{(5.14)}$$
 If $x \notin E$, there exists $\eta_0$ such that $I_E(\eta_0) - \langle x, \eta_0\rangle <
 0$. By taking $\eta = s \eta_0$, with $s \geq 1$, in (5.14), there
 exists a positive constant $M_0$ such that
 $$|\phi(x)| \leq M_0 e^{s(I_E(\eta_0) - \langle x, \eta_0\rangle)}
 e^{s(I_E(\eta_0) - \langle x, \eta_0\rangle)}.$$
 But the second member of this inequality tend to zero when $s$ go
 to the infinity. Then the support of $\phi$ is contained in $E$.\\
 As the support of $\varphi$ is contained in $E$, then from
 Theorem 4.2  we deduce that $$\forall\; y \in \mathbb{R}^d, f(y) =
 \mathcal{F}_D(\varphi)(y).$$
 {\bf Remark 5.1}

 Theorem 5.1 shows that Theorem 4.1 i) which is Theorem 5.2 i) of
 [15],
 can be proved without using
 Theorem 4.4 of [10], chap.3 p. 58.\\
 {\bf Corollary 5.1.} Let $E$ be a $W$-invariant compact convex set
 of $\mathbb{R}^d$. Then for all $f$ in
 $\mathcal{D}(\mathbb{R}^d)$ we have
 $$Supp f \subset E \Longleftrightarrow {}^tV_k(f) \subset
 E. \eqno{(5.15)}$$
 {\bf Proof}

 We deduce the result from the relation
 $$\mathcal{F}_D(f) = \mathcal{F}\;o\; {}^tV_k(f),\;\; f
  \in \mathcal{D}(\mathbb{R}^d),$$
 the Paley-Wiener theorem for the classical Fourier transform
 $\mathcal{F}$ (see [1] Theorem 2.6, p. 17) and Theorem 5.1.\\
 {\bf Corollary 5.2.} Let $E$ be a $W$-invariant compact convex set
 of $\mathbb{R}^d$. Then for all $x \in E$, the support of the
 distribution $\eta_x$ given by the relation (3.7) is contained in
 $E$.\\
 {\bf Proof}

 From (3.3) for all $g$ in $\mathcal{D}(\mathbb{R}^d)$ with
 support in the complementary $E^c$ of $E$ and $f$ in $\mathcal{D}
 (\mathbb{R}^d)$ with support in $E$, we have
 $$\int_{\mathbb{R}^d}V^{-1}_k(g)(x)f(x)dx =
  \int_{\mathbb{R}^d}g(y){}^tV^{-1}_k(f)(y)
 \omega_k(y)dy. \eqno{(5.16)}$$
 But from (5.15) we have
 $$Supp f \subset E \Rightarrow Supp{}^tV^{-1}_k (f)\subset E.$$
 Thus
 $$\int_{\mathbb{R}^d}V^{-1}_k(g)(x)f(x)dx = 0.$$
 This relation implies
 $$\forall\; x \in E,\;\;V^{-1}_k(g)(x) = 0.$$
 But from (3.7) we have
 $$\forall\; x \in E, V^{-1}_k(g)(x) = \langle \eta_x, g\rangle.$$
 Thus the support of $\eta_x$ is contained in $E$.\\
 {\bf Remark 5.2}

 Corollary 5.2 is Propostiion 6.3 of [15] p. 30. The proof of this
 Corollary constitutes another proof of this proposition.

 In the following we give an ameliorated version of the proof of
 Proposition 6.3.

 Let $x \in E$ and $\varepsilon \in ]0,1]$. We consider the
 functions $f_x$ and $f^\varepsilon_x$ given by
 $$\forall\; y \in \mathbb{R}^d , f_x(y) = \frac{e^{-i\langle x, y \rangle}}
 {(1+\|y\|^2)^p},\eqno{(5.17)}$$
 and
 $$\forall\; y \in \mathbb{R}^d,\;\; f^\varepsilon_x(y) = f_x(y)
  \mathcal{F}_D(\mathcal{V}_\varepsilon)(y),\eqno{(5.18)}$$
 with $p \in \mathbb{N}$, such that $p > \gamma + \frac{d}{2} +1$,
 and $\mathcal{V}_\varepsilon$ the  function defined by
 $$\forall\; x \in \mathbb{R}^d,\;\; \mathcal{V}_\varepsilon(x)
  = \frac{1}
 {\varepsilon^{2\gamma+d}}\tilde{\mathcal{V}}(\frac{\|x\|}
 {\varepsilon}),\eqno{(5.19)}$$
 where $\mathcal{V}$ is a radical, positive function in
 $\mathcal{D}(\mathbb{R}^d)$, with support in the ball of center
 $0$ and radius 1, satisfying $\int_{\mathbb{R}^d}\mathcal{V}(x) \omega_k(x) dx =
 1$, and $\tilde{\mathcal{V}}$ the function on $[0, + \infty[$
 given by $\mathcal{V}(x) = \tilde{\mathcal{V}}(\|x\|)$.\\
 As the functions $f_x$ and $f^\varepsilon_x$ belong to $(L^1_k \cap
 L^2_k)(\mathbb{R}^d)$, then from Theorem 4.3 ii), the functions
 $F_x$ and $F^\varepsilon_x$ defined by
 $$\forall\; t \in \mathbb{R}^d,\; F_x(t) = \int_{\mathbb{R}^d}f_x(t)K(iy,x)
 \omega_k(y)dy, \eqno{(5.20)}$$
 $$\forall\; t \in \mathbb{R}^d, F^\varepsilon_x(t) = \int_{\mathbb{R}^d}
 f^\varepsilon_x(t)K(iy,x)\omega_k(y)dy, \eqno{(5.21)}$$
 are in $C(\mathbb{R}^d)$ and from the relation (2.14) and the fact
 that there exists $M > 0$ such that
 $$\forall\; y \in \mathbb{R}^d,\; |\mathcal{F}_D(\mathcal
 {V}_\varepsilon)(y) - 1| \leq \varepsilon M\|y\|^2, \eqno{(5.22)}$$
 we deduce that for all $t \in \mathbb{R}^d$ :
 $$\lim_{\varepsilon \rightarrow 0}F_x^\varepsilon(t) =
  F_x(t).\eqno{(5.23)}$$
 The relation (5.21) can also be written in the form
 $$\forall\; t \in \mathbb{R}^d, F_x^\varepsilon(t) = \int_{\mathbb{R}^d}
 e^{-i\langle x,y\rangle}\mathcal{F}_D(\mathcal{V}_\varepsilon)
 (y) K(iy,t) M_k(y)dy,$$
 where $M_k(y) = \frac{\omega_k(y)}{(1+\|y\|^2)^p}$. Using the
 method applied in the proof of Theorem 5.1, which is the same
 method applied also in the proof of Proposition 6.3,  we deduce
 that
 $$\forall\; t \in \mathbb{R}^d, |F^\varepsilon_x(t)| \leq e^{(I_{E+B_
 \varepsilon}(\eta) - \langle t, \eta\rangle)}
 \int_{\mathbb{R}^d}\frac{|m_k(u + i\eta)|}{(1+\|u\|^2)^p}du,$$
 where $B_\varepsilon$ is the ball of center $0$ and radius
 $\varepsilon$.\\
 As for the proof of Theorem 5.1, we deduce that the support of
 $F^\varepsilon_x$ is contained in $E + B_\varepsilon$.\\ From
 this result and (5.23) we show that
 $$Supp F_x \subset E.$$
 By applying now the remainder of the proof given in [15] p. 32, we
 deduce that Supp$\eta_x \subset E$.
 \subsection{The Dunkl transform of distributions}
 {\bf Notation.} We denote by $\mathcal{H}(\mathbb{C}^d)$ the
 space of entire functions on $\mathbb{C}^d$ which  are slowly
 increasing and of exponential type.
 \vspace{5mm}

 The Dunkl transform of a disgtribution $S$ in
 $\mathcal{E}'(\mathbb{R}^d)$ is defined by
 $$\forall\; y \in \mathbb{R}^d, \mathcal{F}_D(S)(y) =
  \langle S_x, K(-iy,x)\rangle.\eqno{(5.24)}$$
 {\bf Remark 5.3}

 When $S$ is given by a function $g$ in
 $\mathcal{D}(\mathbb{R}^d)$, and denoted by $T_g$, the relation
 (5.24) coincides with (4.1), because $T_g$ is defined
 by
 $$\langle T_g, \varphi\rangle = \int_{\mathbb{R}^d}g(x)
 \varphi(x) \omega_k(x)dx, \varphi \in \mathcal{E}
 (\mathbb{R}^d).\eqno{(5.25)}$$
 The following theorem is proved in [15].\\
 {\bf Theorem 5.2.} The transform $\mathcal{F}_D$ is a topological
 isomorphism from $\mathcal{E}'(\mathbb{R}^d)$ onto
 $\mathcal{H}(\mathbb{C}^d)$.

 We give now the geometrical form of the Paley-Wiener theorem for
 distributions.\\
  {\bf Theorem 5.3.} Let $E$ be a $W$-invariant compact convex set
 of $\mathbb{R}^d$ and $f$ an entire function on $\mathbb{C}^d$. Then
 $f$ is the Dunkl transform of a distribution in
 $\mathcal{E}'(\mathbb{R}^d)$ with support in $E$ if and only if
 there exist a positive constant $C$ and $N \in \mathbb{N}$ such
 that
 $$\forall\; z \in \mathbb{C}^d, |f(z)| \leq C (1 + \|z\|^2)^N
  e^{I_E(Imz)}.\eqno{(5.26)}$$
 {\bf Proof}

 - Necessity condition

 We consider a distribution $S$ in $\mathcal{E}'(\mathbb{R}^d)$
 with support in $E$.

 Let $\mathcal{X}$ be in $\mathcal{D}(\mathbb{R}^d)$ equal to 1 in
 a neighbourhood of $E$, and $\theta$ in $\mathcal{E}(\mathbb{R})$
  such that
  $$\theta(t) = \left\{ \begin{array}{ll}
  1, &\mbox{ if } t \leq 1,\\
  0, &\mbox{ if } t > 2.
  \end{array}\right.$$
  We put $\eta = Im z, z \in \mathbb{C}^d$.\\
  We denote by $\psi_z$ the function defined on $\mathbb{R}^d$ by
  $$\psi_z(x) = \chi(x) K(-ix,z) \theta(\langle x, \eta\rangle - I_E(\eta)).$$
  This function belongs to $\mathcal{D}(\mathbb{R}^d)$ and it is
  equal to $K(-ix,z)$ in a neighbourhood of $E$. Thus
  $$\forall\; z \in \mathbb{C}^d,\;\; \mathcal{F}_D(S)(z)
   = \langle S_x, \psi_z(x)\rangle.\eqno{(5.27)}$$
  As $S$ is with compact support, then it is of finite order $N$.
  Then there exists a positive constant $C_0$ such that
  $$\forall\; z \in \mathbb{C}^d, |\mathcal{F}_D(S)(z)| \leq C_0 \sum_{|p| \leq N}
  \sup_{x \in \mathbb{R}^d}|D^p \psi_z(x)|.\eqno{(5.28)}$$
  Using Leibniz rule, we obtain
  $$\forall\; x \in \mathbb{R}^d, D^p \psi_z(x) = \sum_{q+r+s=p}
  \frac{p!}{q!r!s!} D^q \mathcal{X}(x)D^r K(-ix,z)
  D^s \theta(\langle x,\eta\rangle - I_E(\eta)).\eqno{(5.29)}$$
  We have
  $$\forall\;x \in \mathbb{R}^d, |D^q \mathcal{X}(x)| \leq \mbox{ const. },$$
  and $$\forall\; x \in \mathbb{R}^d, |D^s \theta(\langle x,
  \eta\rangle - I_E(\eta))|
  \leq \mbox{ const. } \|\eta\|^{|s|}.$$
  On the other hand from (2.11) and the fact that $E$ is
  $W$-invariant we have
  $$\forall\; x \in \mathbb{R}^d, |D^r K(-ix,z)| \leq \|z\|^{|r|}
  e^{\langle x, \eta \rangle}.$$
  Using these inequaties and (5.29) we deduce that there exists a
  positive constant $C_1$ such that
  $$\forall\; x \in \mathbb{R}^d, |D^p \psi_z(x)| \leq C_1(1 + \|z\|^2)^N
  e^{\langle x, \eta\rangle}.$$
  From this relation and (5.28) we obtain
  $$\forall\; z \in \mathbb{C}^d, \; |\mathcal{F}_D(S)(z)|
   \leq C_2(1+\|z\|^2)^N \sup
  e^{\langle x, \eta \rangle}, \eqno{(5.30)}$$
  where $C_2$ is a positive constant, and the supremum is
  calculated for
  $$\langle x, \eta\rangle \leq I_E(\eta) + 2.$$
  But this inequality implies
  $$\sup e^{\langle x, \eta\rangle} \leq e^2 e^{I_E(\eta)}.\eqno{(5.31)}$$
  From (5.30) and (5.31) we deduce that there exists a positive
  constant $C$ such that
  $$\forall\; z \in \mathbb{C}^d, |\mathcal{F}_D(S)(z)| \leq C(1 + \|z\|^2)^N
  e^{I_E(\eta)}.\eqno{(5.32)}$$

  - Sufficiency condition

  Let $f$ be an entire function on $\mathbb{C}^d$ satisfiying the
  condition (5.26). We consider the functions $g$ and $g_\varepsilon, \varepsilon \in
   ]0,1]$, given by
$$\forall\; y \in \mathbb{R}^d, g(y) = \frac{f(y)}{(1 +
\|y\|^2)^{p+N}}, \eqno{(5.33)} $$ and $$\forall\; y \in
\mathbb{R}^2, \; g_\varepsilon(y) = g(y) \mathcal{F}_D (
\mathcal{V}_\varepsilon)(y),\eqno{(5.34)}$$ with $p \in
\mathbb{N}$, such that $p > \gamma + \frac{d}{2} +1$, and
$\mathcal{V}_\varepsilon$ the function defined by (5.19).\\ As the
functions $g$ and $g_\varepsilon$ belong to $(L^1_k \cap
L^2_k)(\mathbb{R}^d)$, then from Theorem 4.3 ii) the functions $G$
and $G_\varepsilon$ defined by $$\forall\; x \in \mathbb{R}^d,
G(x) = \int_{\mathbb{R}^d}g(y) K(iy, x)\omega_k(y)dy,$$
$$\forall\; x \in \mathbb{R}^d, G_\varepsilon(x) =
\int_{\mathbb{R}^d}g_\varepsilon(y)K(iy, x)\omega_k(y)dy,$$ are in
$C(\mathbb{R}^d)$ and from (2.14) and (5.22), for all $x \in
\mathbb{R}^d$, we obtain $$\lim_{\varepsilon \rightarrow
0}G_\varepsilon(x) = G(x).\eqno{(5.35)}$$ By making the same proof
as for the Remark 5.1, we deduce that the support of
$G_\varepsilon$ is contained in $E + B_\varepsilon$, where
$B_\varepsilon$ is the ball of center $0$ and radius
$\varepsilon$. From (5.35) we deduce that $$\mbox{ Supp }G \subset
E. $$ We put $$S = (I - \Delta_k)^{p+N})G.$$ where $\Delta_k$ is
the Dunkl Laplacian.\\ Then $S$ is a distribution in
$\mathcal{E}'(\mathbb{R}^d)$ with support in $E$. Moreover, we
have $$\forall\; y \in \mathbb{R}^d, \mathcal{F}_D(S)(y) =
\int_{\mathbb{R}^d}(1+\|y\|^2)^{p+N}G(x)K(-iy, x)\omega_k(x)dx.$$
Using Theorem 4.3 ii) and (5.33) we obtain $$\forall\; y \in
\mathbb{R}^d, \mathcal{F}_D(S)(y) = f(y).$$ This completes the
proof of the Theorem.
\section{Inversion formulas for the Dunkl intertwining operator and for its dual}
 \subsection{The pseudo-differential operators $P$ and $Q$}

 We consider the pseudo-differential operators $P$ and $Q$ defined
 on $\mathcal{D}(\mathbb{R}^d)$ by
 $$\forall\; x \in \mathbb{R}^d,\; P(f)(x) = \frac{2^{2\gamma}}
 {\pi^d c^2_k}\mathcal{F}^{-1}[\omega_k \mathcal{F}(f)](x).\eqno{(6.1)}$$
 $$\forall\;x \in \mathbb{R}^d,\; Q(f)(x) = \frac{2^{2\gamma}}
 {\pi^d c^2_k}\mathcal{F}^{-1}_D[\omega_k \mathcal{F}_D(f)](x).\eqno{(6.2)}$$
 {\bf Proposition 6.1.} For all $f$ in $\mathcal{D}(\mathbb{R}^d)$
 the function $P(f)$ and $Q(f)$ are of class $C^\infty$ on
 $\mathbb{R}^d$ and we have
 $$\frac{\partial}{\partial x_j} P(f)(x) = P(\frac{\partial}{\partial \xi_j}f)(x), j = 1, 2,\cdots,d,\eqno{(6.3)}$$
 $$T_j Q(f)(x) = Q(T_jf)(x),\; j = 1, 2, \cdots,d,\eqno{(6.4)}$$
 {\bf Proof}

 We deduce the results by derivation under the integral sign,
 by using (2.7), (2.8) and the relations
 $$\forall\; y \in \mathbb{R}^d,\; iy_j \mathcal{F}(f)(y) = \mathcal{F}(\frac{\partial}
 {\partial \xi_j}f)(y), \eqno{(6.5)}$$
 $$\forall\; y \in \mathbb{R}^d,\; iy_j \mathcal{F}_D(f)(y) = \mathcal{F}_D(T_jf)(y).\eqno{(6.6)}$$
{\bf Proposition 6.2.} Let $E$ be a $W$-invariant compact convex
set of $\mathbb{R}^d$. Then for all $f$ in
$\mathcal{D}(\mathbb{R}^d)$ we have $$Supp f \subset E \Rightarrow
Supp P(f) \subset E \mbox{ and } Supp Q(f) \subset E. \eqno{(6.7)}
$$ {\bf Proof}

We obtain (6.7) by using the Paley-Wiener theorem for the Fourier
transform $\mathcal{F}$ (see [1] Theorem 2.6, p. 17), Theorem 5.1,
and by applying the methode used in the proof of Theorem 5.1.\\
{\bf Proposition 6.3.} We suppose that $k(\alpha) \in \mathbb{N}$
for all $\alpha \in R_+\;\;\;$. Then for all $f$ in
$\mathcal{D}(\mathbb{R}^d)$ we have $$P(f) = \left[ \prod_{\alpha
\in R_+ }(-1)^{k(\alpha)}\left(\alpha_1 \frac{\partial}{\partial
\xi_1 } +\cdots+ \alpha_d \frac{\partial}{\partial \xi_d
}\right)^{2k(\alpha)}\right]f, \eqno{(6.8)}$$ $$Q(f) =
\left[\prod_{\alpha \in R_+}(-1)^{k(\alpha)}(\alpha_1 T_1 +\cdots+
\alpha_d T_d)^{2k(\alpha)} \right]f.\eqno{(6.9)}$$ {\bf Proof}

As $k(\alpha) \in \mathbb{N}$ for all $\alpha \in R_+$, then for
all $f \in \mathcal{D}(\mathbb{R}^d)$, we have $$\forall \; y \in
\mathbb{R}^d, \omega_k(y) \mathcal{F}(f)(y) = \prod_{\alpha \in
R_+}(\langle \alpha,y
\rangle)^{2k(\alpha)}\mathcal{F}(f)(y),\eqno{(6.10)} $$ and
$$\forall\; y \in \mathbb{R}^d,\; \omega_k(y)\mathcal{F}_D(f)(y) =
\prod_{\alpha \in R_+}(\langle
\alpha,y\rangle)^{2k(\alpha)}\mathcal{F}_D(f)(y). \eqno{(6.11)} $$
But $$\forall\; y \in \mathbb{R}^d, \langle \alpha,y\rangle
\mathcal{F}(f)(y) = \mathcal{F}\left[- i \left(\alpha_1
\frac{\partial}{\partial \xi_1} + \cdots+ \alpha_d
\frac{\partial}{\partial \xi_d} \right)f\right](y),
\eqno{(6.12)}$$ and using (2.7), (2.8) we deduce that $$\forall\;
y \in \mathbb{R}^d, \langle \alpha, y\rangle \mathcal{F}_D(f)(y) =
\mathcal{F}_D[-i(\alpha_1 T_1 +\cdots+ \alpha_d
T_d)f](y).\eqno{(6.13)} $$ From (6.10), (6.11) and (6.12), (6.13)
we obtain $$\forall\; y \in \mathbb{R}^d, \omega_k (y)
\mathcal{F}(f)(y) = \mathcal{F}\left[\prod_{\alpha \in
R_+}(-1)^{k(\alpha)}(\alpha_1 \frac{\partial}{\partial \xi_1} +
\cdots+ \alpha_d \frac{\partial}{\partial \xi_d})^{2k(\alpha)}f
\right](y),$$ and $$\forall\; y \in \mathbb{R}^d, \omega_k(y)
\mathcal{F}_D(f)(y) = \mathcal{F}_D\left[\prod_{\alpha \in R_+ }
(-1)^{k(\alpha)}(\alpha_1 T_1 +\cdots+ \alpha_dT_d)^{2k(\alpha)}f
\right](y).$$ These relations, the inversion formula for the
Fourier transform $\mathcal{F}$, and Theorem 4.2 imply (6.8),
(6.9).\\ {\bf Theorem 6.1.} For all $f$ in
$\mathcal{D}(\mathbb{R}^d)$ we have the following  transmutation
relation $$\forall\; x \in \mathbb{R}^d,\; P({}^tV_k(f))(x) =
{}^tV_k(Q(f))(x).\eqno{(6.14)}$$ {\bf Proof}

From Proposition 6.2 and the relations (6.1), (4.3) we have
\begin{eqnarray*}
\forall\; x \in \mathbb{R}^d, P({}^tV_k(f))(x) &=&
\frac{2^{2\gamma}}{\pi^d c^2_k}\mathcal{F}^{-1}\Big[\omega_k
\mathcal{F}\,o\, {}^tV_k(f)\Big](x),\\ &=&
\frac{2^{2\gamma}}{\pi^d c^2_k}{}^tV_k \Big\{
\mathcal{F}_D[\omega_k \mathcal{F}_D(f)]\Big\}(x).
\end{eqnarray*}
Then (6.2) implies $$\forall\; x \in \mathbb{R}^d,\,
P({}^tV_k(f))(x) = {}^tV(Q(f))(x).$$
\subsection{Inversion formulas for the Dunkl intertwining operator and for its dual}

In this subsection we give inversion formulas for the operators
$V_k$ and ${}^tV_k$ and we deduce the expressions of the
representing distributions of the operators $V^{-1}_k$ and
${}^tV_k^{-1}$.\\ {\bf Theorem 6.2.} For all $f$ in
$\mathcal{D}(\mathbb{R}^d)$ we have $$\forall\;x \in
\mathbb{R}^d,\; {}^tV^{-1}_k(f)(x) = V_k(P(f))(x).\eqno{(6.15)}$$
{\bf Proof}

From Theorem 3.2 ii),  for $f$ in $\mathcal{D}(\mathbb{R}^d)$ the
function ${}^tV^{-1}_k(f)$ belongs to $\mathcal{D}(\mathbb{R}^d)$.
Then from Theorem 4.2 we have $$\forall\; x \in \mathbb{R}^d,
{}^tV^{-1}_k(f)(x) = \frac{c^2_k}{2^{2\gamma
+d}}\int_{\mathbb{R}^d}K(iy, x)
\mathcal{F}_D({}^tV^{-1}_k)(f))(y)\omega_k(y)dy.\eqno{(6.16)}$$
But from the relations (4.3), (3.1), we have $$\forall\; y \in
\mathbb{R}^d,\; \mathcal{F}_D({}^tV^{-1}_k(f))(y) =
\mathcal{F}(f)(y),$$ and $$\forall  y \in \mathbb{R}^d,\; K(iy,x)
= \mathcal{F}(\check{\mu}_x)(y),$$ where $\check{\mu}_x$ is the
probability measure given by $$\int_{\mathbb{R}^d}f(t)
d\check{\mu}_x (t) = \int_{\mathbb{R}^d}f(-t)d\mu_x(t), f \in
C(\mathbb{R}^d).$$ Thus (6.16) can also be written in the form
$$\forall\; x \in \mathbb{R}^d, {}^tV^{-1}_k(f)(x) =
\frac{c^2_k}{2^{2\gamma+d}}
\int_{\mathbb{R}^d}\mathcal{F}(\check{\mu}_x)(y)\omega_k(y)
\mathcal{F}(f)dy.\eqno{(6.17)}$$ But from Propostion 6.2 and (6.1)
we have $$\forall\; y \in \mathbb{R}^d, \omega_k(y)
\mathcal{F}(f)(y) =
\frac{2^{2\gamma}}{\pi^dc^2_k}\mathcal{F}(P(f)))(y).$$ Then by
using (6.17) and the properties of the Fourier transform
$\mathcal{F}$ we obtain
\begin{eqnarray*}
\forall\; x \in \mathbb{R}^d, {}^tV^{-1}_k(f)(x) &=&
\frac{1}{(2\pi)^d}\int_{\mathbb{R}^d} \mathcal{F}(\check{\mu}(y)
\mathcal{F}(P(f))(y)dy,\\ &=&
\frac{1}{(2\pi)^d}\int_{\mathbb{R}^d}\mathcal{F}(\check{\mu}_x
\ast P(f))(y)dy,
\end{eqnarray*}
where $\ast$ is the classical convolution product on
$\mathbb{R}^d$ of a measure and a function.\\ By using Proposition
6.2, the fact that $\check{\mu}_x$ is a probability measure on
$\mathbb{R}^d$ and the inversion formula for the Fourier transform
$\mathcal{F}$,  we deduce that $$\forall\; x \in \mathbb{R}^d,
{}^tV^{-1}_k(f)(x) = \check{\mu}_x \ast P(f)(0).$$ But from (3.1)
we have $$\check{\mu}_x \ast P(f)(0) = \int_{\mathbb{R}^d}P(f)(-t)
d\check{\mu}_x(t) = V_k(P(f))(x). $$ Thus $$\forall\; x \in
\mathbb{R}^d, {}^tV^{-1}_k(f)(x) = V_k(P(f))(x).$$ {\bf Corollary
6.1.} For all $f$ in $\mathcal{D}(\mathbb{R}^d)$ we have the
following relation $$\forall\; x \in \mathbb{R}^d, V_k(P^2(f))(x)
= QV_k(P(f))(x).\eqno{(6.18)} $$ {\bf Proof}

We deduce (6.18) from Theorem 6.2, 6.3. \\ {\bf Theorem 6.4.} For
all $f$ in $\mathcal{D}(\mathbb{R}^d)$ we have $$\forall\;  x \in
\mathbb{R}^d,\; V^{-1}_k(f)(x) = {}^tV_k(Q(f))(x).\eqno{(6.19)} $$
{\bf Proof}

We obtain the result by using Proposition 6.3 and Theorems 6.1,
6.2.\\ {\bf Notation}

We denote by $\mathcal{D}'(\mathbb{R}^d)$ the space of
distributions on $\mathbb{R}^d$. It is the topological dual of
$\mathcal{D}(\mathbb{R}^d)$.\\ {\bf Definition 6.1.} We define the
transposed operators ${}^tP$ and ${}^tQ$ of the operators $P$ and
$Q$ on $\mathcal{D}'(\mathbb{R}^d)$ by $$\langle {}^tP(S),f\rangle
= \langle S, P(f) \rangle , f \in \mathcal{D}(\mathbb{R}^d),
\eqno{(6.20)} $$ $$\langle {}^tQ(S), f\rangle = \langle
S,Q(f)\rangle, f \in \mathcal{D}(\mathbb{R}^d).\eqno{(6.21)}$$
{\bf Proposition 6.4.} We suppose that $k(\alpha) \in \mathbb{N}$
for all $\alpha \in R_+$. Then for all $S \in
\mathcal{D}'(\mathbb{R}^d)$ we have $${}^tP(S) =
\left[\prod_{\alpha \in R_+} (\alpha \frac{\partial}{\partial
\xi_1} + \cdots + \alpha_d \frac{\partial}{\partial
\xi_d})^{2k(\alpha)}\right]S,\eqno{(6.22)}$$ $${}^tQ(S) =
\left[\prod_{\alpha \in R_+}(\alpha T_1 + \cdots+
\alpha_dT_d)^{2k(\alpha)} \right]S,\eqno{(6.23)}$$ where $T_j, j =
1, 2,\cdots,d$, are the Dunkl operators defined on
$\mathcal{D}'(\mathbb{R}^d)$ by $$\langle T_j S,f\rangle = -
\langle S, T_jf\rangle, f\in \mathcal{D}(\mathbb{R}^d).$$ {\bf
Theorem 6.5.} The representing distributions $\eta_x$ and $Z_x$ of
the inverse of the Dunkl intertwning operator and of its dual, are
given by $$\forall\; x \in \mathbb{R}^d, \; \eta_x =
{}^tQ(\nu_x),\eqno{(6.24)}$$ and $$\forall\;x \in \mathbb{R}^d,
Z_x = {}^tP(\mu_x).\eqno{(6.25)}$$ where $\mu_x$ and $\nu_x$ are
the representing measues of the Dunkl intertwining operator $V_k$
and of its dual ${}^tV_k$.\\ {\bf Proof}

- From (3.4), for all $f$ in $\mathcal{D}(\mathbb{R}^d)$ we have
$$\forall\; x \in \mathbb{R}^d,\; {}^tV_k(Q(f))(x) = \langle\nu_x,
Q(f)\rangle = \langle {}^tQ(\nu_x),f\rangle.\eqno{(6.26)} $$ On
the other hand from (3.7) : $$\forall\; x \in \mathbb{R}^d,
V^{-1}_k(f)(x) = \langle \eta_x,f\rangle.$$ We obtain (6.24) from
this relation, (6.26) and (6.19).

- By using (3.1), for all $f$ in $\mathcal{D}(\mathbb{R}^d)$ the
relation (6.15) can also be written in the form $$\forall\; x \in
\mathbb{R}^d, {}^tV^{-1}_k(f)(x) = \langle \mu_x, P(f)\rangle =
\langle{}^tP(\mu_x),f\rangle.\eqno{(6.27)} $$ But from (3.10) we
have $$\forall\; x \in \mathbb{R}^d, {}^tV^{-1}(f)(x) = \langle
Z_x,f\rangle.$$ We deduce (6.25) from this relation and (6.27).\\
{\bf Remark 6.1}

When $k(\alpha) \in \mathbb{N}$ for all $\alpha \in R_+$, we have
$$\forall\; x \in \mathbb{R}^d, \eta_x = \left[\prod_{\alpha \in
R_+ }\left(\alpha_1T_1+\cdots+ \alpha_d T_d
\right)^{2k(\alpha)}\right](\nu_x),\eqno{(6.28)}$$ and $$\forall\;
x \in \mathbb{R}^d, Z_x = \left[\prod_{\alpha \in R_+}(\alpha_1
\frac{\partial}{\partial \xi_1}+\cdots+ \alpha_d
\frac{\partial}{\partial
\xi_d})^{2k(\alpha)}\right](\mu_x).\eqno{(6.29)}$$ {\bf
R\'ef\'erences}
\begin{enumerate}
\item Chazarain. J. and Piriou. A. (1982). Introduction to the
theory of linerar partial differential equations. North Holland
Publishing Company - Amsterdam, New-York. Oxford.
\item Van Diejen, J. F. (1997) . Confluent hypergeoemtric
orthogonal polynomials related to the rational quantum Calogero
system with harmonic confinement.Comm. Math. Phys., {\bf 188},
467-497.
\item Dunkl, C. F. (1989) Differential-difference operators
associated to reflection groups. Trans. Amer.Math. Soc, {\bf 311},
167-183
\item Dunkl,C.F. (1991). Integral kernels with reflection group
invariance. Can. J. Math., {\bf 43}, 1213-1227.
\item Dunkl, C.F. (1992). Hankel transforms associated to finite
reflection groups. Contemp. Math., {\bf 138}, 123-138.
\item Heckman, G.J. (1991). An elementary approach to the
hypergeometric shift operators of Opdam. Invent. Math., {\bf 103},
341-350.
\item Humphreys, J.E. (1990). Reflection Groups and Coxeter Groups.
Cambridge Univ. Press., Cambridge, England.
\item Hikami, K. (1996), Dunkl operators formalism for quantum
many-body problems associated with classical root systems. J.
Phys. Soc. Japan, {\bf 65}, 394-401.
\item de Jeu, M.F.E (1993). The Dunkl transform. Invent. Math.,
{\bf 113}, 147-162.
\item de Jeu, M.F.E. (1994).Dunkl operators. Thesis, University of
Amsterdam.
\item Kakei, S. (1996). Common algebraic structure for the
Calogero-Stherland models. J. Phys. A, {\bf 29}, 619-624.
\item Lapointe, M. and Vinet, L. (1996). Exact operator solution
of the Calogero-Surherland model. Comm. Math. Phys., {\bf 178},
425-452.
\item Rosler, M. (1999). Positivity of Dunkl's intertwining
operator. Duke. Math. J., {\bf 98}, 445-463.
\item Trim\`eche, K. (2001). The Dunkl intertwining operator on
spaces of functions and distributions and integral representation
of its dual. Integral Transforms and Special Functions, {\bf 12}
(4), 349-374.
\item Trim\`eche.K. (2002). Paley-Wiener theorems for the Dunkl
transform and Dunkl translation operators. Integral Transform and
Special Functions, {\bf 13}, 17-38.
\item Xu. Y. (1997). Orthogonal polynomials for a family of
product weight function on the sphere. Can. J. Math. 49, 175-192.
\end{enumerate}
 \end{document}